\documentclass[12pt]{article}
\usepackage{latexsym, amssymb, amsmath, amscd, amsfonts, epsfig, graphicx, colordvi,verbatim,ifpdf}
\usepackage{amsfonts, amsmath, amssymb}
\usepackage{amssymb,amsfonts,amsmath,latexsym,epsfig,cite, psfrag,eepic,color}
\usepackage{amscd,graphics}
\usepackage{latexsym, amssymb,  amsmath,amscd, amsfonts, epsfig, graphicx, colordvi,amsthm}

\usepackage{graphicx}
\usepackage{color}
\usepackage{ifpdf}
\usepackage{fancybox}

\usepackage{float}
\usepackage{cases}

\newtheorem{thm}{Theorem}[section]
\newtheorem{pro}[thm]{Proposition}

\newtheorem{defi}[thm]{Definition}

\newtheorem{coro}[thm]{Corollary}

\def\pf{\noindent{\it Proof.} }
\def\qed{\nopagebreak\hfill{\rule{4pt}{7pt}}
\medbreak}

\numberwithin{equation}{section}

\def\qed{\nopagebreak\hfill{\rule{4pt}{7pt}}
\medbreak}

\newlength{\boxedparwidth}
  {\begin{center} \begin{tabular}{|@{\hspace{.315in}}c@{\hspace{.15in}}|}
                  \hline \\ \begin{minipage}[t]{\boxedparwidth}
                  \setlength{\parindent}{.25in}}%
  {\end{minipage} \\ \\ \hline \end{tabular} \end{center}}

\begin{document}
\begin{center}

 {\large \bf Ramanujan-type Congruences for Overpartitions Modulo $5$}

\end{center}

\begin{center}
{  William Y.C. Chen}$^{a,b}$, {Lisa H. Sun}$^{a,*}$, {Rong-Hua Wang}$^{a}$  and
  {Li Zhang}$^{a}$ \vskip 2mm

   $^{a}$Center for Combinatorics, LPMC-TJKLC\\
   Nankai University, Tianjin 300071, P. R. China\\[6pt]
   $^{b}$Center for Applied Mathematics\\
Tianjin University,  Tianjin 300072, P. R. China\\[12pt]

    chen@nankai.edu.cn, $^{*}$sunhui@nankai.edu.cn, \\
    wangwang@mail.nankai.edu.cn, zhangli427@mail.nankai.edu.cn
\end{center}

\vskip 6mm \noindent {\bf Abstract.}  Let $\overline{p}(n)$ denote the number of overpartitions of $n$. Hirschhorn and Sellers showed that $\overline{p}(4n+3)\equiv 0 \pmod{8}$ for   $n\geq 0$. They also conjectured that $\overline{p}(40n+35)\equiv 0 \pmod{40}$ for $n\geq 0$.  Chen and Xia proved this conjecture by using the $(p,k)$-parametrization of theta functions given by Alaca, Alaca and Williams. In this paper, we show that
$\overline{p}(5n)\equiv (-1)^{n}\overline{p}(4\cdot 5n) \pmod{5}$ for $n \geq 0$ and  $\overline{p}(n)\equiv (-1)^{n}\overline{p}(4n)\pmod{8}$ for $n \geq 0$ by using the relation of the generating function of $\overline{p}(5n)$ modulo $5$ found by Treneer and the $2$-adic expansion of the generating function of $\overline{p}(n)$ due to Mahlburg. As a consequence, we deduce that $\overline{p}(4^k(40n+35))\equiv 0 \pmod{40}$ for $n,k\geq 0$.
Furthermore, applying the Hecke operator on $\phi(q)^3$ and the fact that $\phi(q)^3$ is a Hecke eigenform, we obtain an infinite family of congrences $\overline{p}(4^k \cdot5\ell^2n)\equiv 0 \pmod{5}$, where $k\ge 0$ and $\ell$ is a prime such that $\ell\equiv3 \pmod{5}$ and $\left(\frac{-n}{\ell}\right)=-1$.
Moreover, we show that $\overline{p}(5^{2}n)\equiv \overline{p}(5^{4}n) \pmod{5}$ for $n \ge 0$. So we are led to the congruences $\overline{p}\big(4^k5^{2i+3}(5n\pm1)\big)\equiv 0 \pmod{5}$ for $n, k, i\ge 0$.
In this way, we obtain various Ramanujan-type congruences for $\overline{p}(n)$ modulo $5$ such as $\overline{p}(45(3n+1))\equiv 0 \pmod{5}$ and $\overline{p}(125(5n\pm 1))\equiv 0 \pmod{5}$ for $n\geq 0$.



\section{Introduction}

 The objective of this paper is to use  half-integral weight modular forms to derive three infinite families of congruences for overpartitions modulo $5$.

 Recall that a partition of a nonnegative integer $n$ is a nonincreasing sequence of positive integers whose sum is $n$. An overpartition of $n$ is a partition of $n$ where the first occurrence of each distinct part may be overlined. We denote the number of overpartitions of $n$ by $\overline{p}(n)$. We set $\overline{p}(0)=1$ and $\overline{p}(n)=0$ if $n<0$. For example, there are eight  overpartitions of $3$
\[3,\ \bar{3},\ 2+1,\ \bar{2}+1,\ 2+\bar{1},\ \bar{2}+\bar{1},\ 1+1+1,\ \bar{1}+1+1.\]

Overpartitions arise in combinatorics \cite{Corteel-Lovejoy-2003}, $q$-series \cite{Corteel-Hitczenko-2004}, symmetric functions \cite{Berndt-2006}, representation theory \cite{Kang-Kwon-2004},
mathematical physics \cite{1-Fortin-Jacob-Mathieu-2005,Fortin-Jacob-Mathieu-2005} and number theory \cite{Lovejoy-2004,Lovejoy-Mallet-2011}.
They are also called standard MacMahon diagrams, joint partitions,  jagged partitions or dotted partitions.

 Corteel and Lovejoy \cite{Corteel-Lovejoy-2003} showed that the generating function of $\overline{p}(n)$ is given by
\[
\sum_{n\geq0}\overline{p}(n)q^n= \frac{(-q;q)_\infty}{(q;q)_\infty}.
\]
Recall that the generating function of $\overline{p}(n)$ can be expressed as
\begin{equation*}
\sum_{n\geq0}\overline{p}(n)q^n=\frac{1}{\phi(-q)},
\end{equation*}
where $\phi(q)$ is Ramanujan's  theta function as defined by
\begin{equation}\label{sum}
\phi(q)=\sum_{n=-\infty}^{\infty}q^{n^2},
\end{equation}
see Berndt \cite{Berndt-2006}.

On the other  hand,  the generating function of $\overline{p}(n)$ has the following $2$-adic expansion
\begin{align}\label{L-2-1}
\sum_{n\geq0}\overline{p}(n)q^n=1+\sum_{k=1}^{\infty}2^k\sum_{n=1}^{\infty}(-1)^{n+k}c_k(n)q^n,
\end{align}
 where $c_k(n)$ denotes the number of representations of $n$
 as a sum of $k$ squares of positive integers. The above
 $2$-adic expansion \eqref{L-2-1}  is useful
 to derive congruences for $\overline{p}(n)$ modulo powers of $2$, see, for example \cite{Kim-2008-1,Kim-2009,Mahlburg-2004}.

By employing  dissection
formulas, Fortin, Jacob and Mathieu \cite{Fortin-Jacob-Mathieu-2005}, Hirschhorn and Sellers \cite{1-Hirschhorn-Sellers-2005} independently derived various Ramanujan-type congruences for $\overline{p}(n)$, such as
\begin{equation}\label{4n+3}
\overline{p}(4n+3)\equiv 0 \pmod{8}.
\end{equation}
Hirschhorn and Sellers \cite{1-Hirschhorn-Sellers-2005} proposed the following conjectures
\begin{align}
&\overline{p}(27n+18)\equiv 0 \pmod{12},  \label{1-r-1}\\[5pt]
&\overline{p}(40n+35)\equiv 0 \pmod{40}.   \label{1-r-2}
\end{align}
They also conjectured that if $\ell$ is prime and $r$ is a quadratic nonresidue modulo $\ell$ then
\begin{align}\label{1-r-3}
\overline{p}(\ell n+r)\equiv\left\{
 \begin{array}{ll}
   0 \pmod{8} &\mbox{if $\ell\equiv\pm1 \pmod{8}$,}\\[5pt]
   0 \pmod{4} &\mbox{if $\ell\equiv\pm3 \pmod{8}$.}
     \end{array}
   \right.
\end{align}

By using the $3$-dissection formula for $\phi(-q)$, Hirschhorn and Sellers \cite{2-Hirschhorn-Sellers-2005} proved \eqref{1-r-1} and obtained  a family of congruences
\begin{align*}
\overline{p}(9^\alpha(27n+18))\equiv 0 \pmod{12},
\end{align*}
where $n,\alpha \geq 0$.

Employing the $2$-dissection formulas of theta functions due to Ramanujan, Hirschhorn and Sellers \cite{1-Hirschhorn-Sellers-2005}, Chen and Xia \cite{Chen-Xia-2013} obtained a generating function of $\overline{p}(40n+35)$ modulo $5$. Using the $(p,k)$-parametrization of theta functions given by Alaca, Alaca and Williams \cite{AAW-2006}, they showed that
\begin{align}\label{40n+35-mod5}
\overline{p}(40n+35)\equiv 0 \pmod{5}.
\end{align}
This proves Hirschhorn and Sellers' conjecture \eqref{1-r-2} by combining congruence \eqref{4n+3}.
Applying the $2$-adic expansion \eqref{L-2-1}, Kim \cite{Kim-2009} proved  \eqref{1-r-3} and obtained  congruence properties  of $\overline{p}(n)$ modulo 8.

For powers of $2$, Mahlburg \cite{Mahlburg-2004} showed that $\overline{p}(n)\equiv0 \pmod{64}$ holds for a set of integers of arithmetic density $1$. Kim \cite{Kim-2008-1} showed that $\overline{p}(n)\equiv0 \pmod{128}$ holds for a set of integers of arithmetic density $1$. For the modulus $3$,
 by using the fact that $\phi(q)^5$ is a Hecke eigenform in the half-integral weight modular form space $M_{\frac{5}{2}}(\tilde{\Gamma}_0(4))$, Lovejoy and Osburn \cite{LOQ11} proved that
\[
\overline{p}(3\ell^3 n)\equiv0 \pmod{3},
\]
where $\ell\equiv2 \pmod{3}$ is an odd prime and $\ell \nmid n$.
Moreover, by utilizing  half-integral weight modular forms, Treneer \cite{Treneer-2006} showed that for a prime $\ell$ such that $\ell\equiv-1 \pmod{5}$, 
\begin{align*}
\overline{p}(5 \ell^3 n)\equiv0 \pmod{5},
\end{align*}
for all $n$ coprime to $\ell$.

In this paper, we   establish the following two congruence relations for overpartitions modulo $5$ and modulo $8$ by using a relation of the generating function of $\overline{p}(5n)$ modulo $5$  and applying the $2$-adic expansion \eqref{L-2-1}.

\begin{thm}\label{rel-mod5}
For   $n\geq0$, we have
\begin{align}\label{recmod5}
&\overline{p}(5n)\equiv (-1)^{n}\overline{p}(4\cdot 5n) \pmod{5}.
\end{align}
\end{thm}

\begin{thm}\label{rel-mod8}
For $n\geq0$, we have
\begin{align}\label{formula-mod8}
\overline{p}(n)\equiv (-1)^{n}\overline{p}(4n)\pmod{8}.
\end{align}
\end{thm}

Combining the above two congruence relations with congruences \eqref{4n+3} and \eqref{40n+35-mod5}, we arrive at a family of congruences modulo 40.

\begin{coro}\label{gen-HS}
For  $n,k\geq0$, we have
\begin{align}\label{m1}
\overline{p}(4^k(40n+35))\equiv 0 \pmod{40}.
\end{align}
\end{coro}

 Based on the Hecke operator on $\phi(q)^3$ and the fact that $\phi(q)^3$ is a Hecke eigenform in   $M_{\frac{3}{2}}(\tilde{\Gamma}_0(4))$, we obtain
  a family of congruences for  overpartitions modulo $5$.

\begin{thm}\label{thm-2}  Let $\left(\frac{\cdot}{\ell}\right)$ denote the Legendre symbol.
Assume that  $k$ is a nonnegative integer and $\ell$ is a prime with $\ell\equiv3 \pmod{5}$. Then we have
\begin{align*}
\overline{p}(4^k \cdot5\ell^2n)\equiv 0 \pmod{5},
\end{align*}
where $n$ is a nonnegative integer such that $\left(\frac{-n}{\ell}\right)=-1$.
\end{thm}

 Using the properties of  the Hecke operator $T_{\frac{3}{2},16}(\ell^2)$ and the Hecke eigenform $\phi(q)^3$, we are led to another congruence relation for overpartitions modulo $5$.

\begin{thm}\label{relation-mod5}
For  $n\geq0$, we have
\begin{align}\label{T-5-1}
&\overline{p}(5^{2}n)\equiv \overline{p}(5^{4}n) \pmod{5}.
\end{align}
\end{thm}

Combining \eqref{recmod5} and \eqref{T-5-1}, we find the
 following family of congruences modulo 5.

\begin{coro}\label{coro-2}
For $n,k,i\geq0$, we have
\begin{align}\label{infi-con-3}
\overline{p}\big(4^k5^{2i+3}(5n\pm1)\big)\equiv 0 \pmod{5}.
\end{align}
\end{coro}

\section{Preliminaries}

To make this paper self-contained, we recall some definitions and notation on half-integral weight modular forms. For more details, see \cite{Chen-Du-Hou-Sun-2012,Treneer-2006,Koblitz-1993,
Ono-2004,Shimura-1973}.

Let $k$ be an odd positive integer and $N$ be a positive integer with
$4|N$. We use $M_{\frac{k}{2}}(\tilde{\Gamma}_0(N))$ to denote
the space of holomorphic modular forms on $\Gamma_0(N)$ of weight $\frac{k}{2}$.

\begin{defi}
Let
\[
f(z)=\sum_{n\geq0}a(n)q^n
\]
be a modular form in $M_\frac{k}{2}(\tilde{\Gamma}_0(N))$.
For any odd prime $\ell\nmid N$, the action of
the Hecke operator $T_{\frac{k}{2},N}(\ell^2)$  on $f(z)\in M_{\frac{k}{2}}(\tilde{\Gamma}_0(N))$ is given by
\begin{align} \label{hec-1}
f(z)\mid T_{\frac{k}{2},N}(\ell^2)=\sum_{n\ge 0}\left(a(\ell^2n)+\Big(\frac{(-1)^{\frac{k-1}{2}n}}{\ell}\Big)
\ell^{\frac{k-3}{2}}a(n)+
\ell^{k-2}a\Big(\frac{n}{\ell^2}\Big)\right)q^n,
\end{align}
where $a(\frac{n}{\ell^2})=0$ if $n$ is not divisible by $\ell^2$.
\end{defi}

The following proposition says that the Hecke operator  $T_{\frac{k}{2},N}(\ell^2)$ maps
the modular form space $M_{\frac{k}{2}}(\tilde{\Gamma}_0(N))$ into itself.

\begin{pro}
Let $\ell$ be an odd prime and $f(z)\in M_{\frac{k}{2}}(\tilde{\Gamma}_0(N))$, then
\[
f(z)\mid T_{\frac{k}{2},N}(\ell^2) \in M_{\frac{k}{2}}(\tilde{\Gamma}_0(N)).
\]
\end{pro}

A Hecke eigenform associated with the Hecke operator $T_{\frac{k}{2},N}(\ell^2)$ is defined as follows.
\begin{defi}
A half-integral weight modular form $f(z)\in M_{\frac{k}{2}}(\tilde{\Gamma}_0(4N))$ is called a Hecke eigenform for the Hecke operator $T_{\frac{k}{2},N}(\ell^2)$, if for every prime $\ell \nmid 4N$ there exists a complex number $\lambda(\ell)$ for which
\[
f(z)\mid T_{\frac{k}{2},N}(\ell^2)=\lambda(\ell) f(z).
\]
\end{defi}

For the space of half-integral weight modular forms on $\Gamma_0(4)$, we have the following dimension formula.

\begin{pro}
We have
\[
\dim M_{\frac{k}{2}}(\tilde{\Gamma}_0(4))=1+\Big\lfloor \frac{k}{4} \Big\rfloor.
\]
\end{pro}

By the above  dimension formula, we see that $\dim M_{\frac{3}{2}}(\tilde{\Gamma}_0(4))=1$. From the fact that $\phi(q)^3 \in M_{\frac{3}{2}}(\tilde{\Gamma}_0(4))$, it is easy to deduce that
\begin{align}\label{eigenform-1}
\phi(q)^3\mid T_{\frac{3}{2},4}(\ell^2)=(\ell+1)\phi(q)^3,
\end{align}
 see, for example \cite[P. 18]{Treneer-2006}.

\section{Proofs of Theorem \ref{rel-mod5} and Theorem \ref{rel-mod8}}

In this section, we give proofs of Theorem \ref{rel-mod5} and Theorem \ref{rel-mod8} by using a relation of the generating function of $\overline{p}(5n)$ modulo $5$ and the $2$-adic expansion \eqref{L-2-1} of $\overline{p}(n)$.

{\noindent \emph{Proof of Theorem \ref{rel-mod5}}}. Recall  the
following $2$-dissection formula for $\phi(q)$,
\begin{align}\label{P-L-3-2}
\phi(q)=\phi(q^4)+2q\psi(q^8),
\end{align}
where
\begin{equation*}
\psi(q)=\sum_{n=0}^{\infty}q^{\frac{n^2+n}{2}},
\end{equation*}
see, for example, Hirschhorn and Sellers \cite{1-Hirschhorn-Sellers-2005}.
Replacing $q$ by $-q$,  \eqref{P-L-3-2} becomes
\begin{align}
\phi(-q)=\phi(q^4)-2q\psi(q^8).\label{dis-q}
\end{align}
We now consider the generating function of $\overline{p}(5n)$ modulo 5.
The following relation is due to Treneer \cite[p. 18]{Treneer-2006},
\begin{equation}\label{L-1-1}
\sum_{n\geq0}\overline{p}(5 n)q^n \equiv \phi(-q)^3\pmod{5}.
\end{equation}
Plugging \eqref{dis-q} into \eqref{L-1-1} yields that
\begin{align}\label{P-L-3-3}
\sum_{n\geq0}\overline{p}(5 n)q^n
\equiv\phi(q^4)^3-q\phi(q^4)^2\psi(q^8)+2q^2\phi(q^4)\psi(q^8)^2-3q^3\psi(q^8)^3 \pmod{5}.
\end{align}
Extracting the terms of $q^{4n+i}$ for $i=0,1,2,3$  on both sides of \eqref{P-L-3-3} and   setting $q^4$ to $q$,  we obtain
\begin{align}
&\sum_{n\geq0}\overline{p}(20n)q^n\equiv\phi(q)^3 \pmod{5},\label{P-L-3-4-1}\\
&\sum_{n\geq0}\overline{p}(20n+5)q^n\equiv-\phi(q)^2\psi(q^2) \pmod{5},\label{P-L-3-4-2}\\
&\sum_{n\geq0}\overline{p}(20n+10)q^n\equiv2\phi(q)\psi(q^2)^2 \pmod{5},\label{P-L-3-4-3}\\
&\sum_{n\geq0}\overline{p}(20n+15)q^n\equiv-3\psi(q^2)^3 \pmod{5}.\label{P-L-3-4-4}
\end{align}
Substituting the $2$-dissection formula \eqref{P-L-3-2} 
 into \eqref{P-L-3-4-1}, we find that
\begin{align}\label{eq-p20}
\sum_{n\geq0}\overline{p}(20 n)q^n \equiv\phi(q^4)^3+q\phi(q^4)^2\psi(q^8)+2q^2\phi(q^4)\psi(q^8)^2+3q^3\psi(q^8)^3 \pmod{5}.
\end{align}
Extracting the terms of $q^{4n+i}$ for $i=0,1,2,3$  on both sides of \eqref{eq-p20} and setting $q^4$ to $q$,  we obtain
\begin{align}
&\sum_{n\geq0}\overline{p}(4\cdot20n)\equiv\phi(q)^3 \pmod{5},\label{P-L-3-7-1}\\[5pt]
&\sum_{n\geq0}\overline{p}(4\cdot(20n+5))\equiv\phi(q)^2\psi(q^2) \pmod{5},\label{P-L-3-7-2}\\[5pt]
&\sum_{n\geq0}\overline{p}(4\cdot(20n+10))\equiv2\phi(q)\psi(q^2)^2 \pmod{5},\label{P-L-3-7-3}\\[5pt]
&\sum_{n\geq0}\overline{p}(4\cdot(20n+15))\equiv3\psi(q^2)^3 \pmod{5}.\label{P-L-3-7-4}
\end{align}
Comparing the equations \eqref{P-L-3-4-1}--\eqref{P-L-3-4-4} with \eqref{P-L-3-7-1}--\eqref{P-L-3-7-4}, we deduce that
\begin{align*}
&\overline{p}(5\cdot(4n))\equiv\overline{p}(4\cdot 5\cdot4n) \pmod{5},\\[5pt]
&\overline{p}(5\cdot(4n+1))\equiv-\overline{p}(4\cdot 5\cdot(4n+1)) \pmod{5},\\[5pt]
&\overline{p}(5\cdot(4n+2))\equiv\overline{p}(4\cdot 5\cdot(4n+2)) \pmod{5},\\[5pt]
&\overline{p}(5\cdot(4n+3))\equiv-\overline{p}(4\cdot 5\cdot(4n+3)) \pmod{5}.
\end{align*}
So we conclude that
\begin{align*}
\overline{p}(5n)\equiv(-1)^n\overline{p}(4\cdot 5n) \pmod{5}.
\end{align*}
 This completes the proof. \qed

We note that extracting the terms of odd powers of $q$ on both sides of \eqref{P-L-3-4-4} leads to the congruence
 $\overline{p}(40n+35)\equiv 0 \pmod{5}$ due to Chen and Xia \cite{Chen-Xia-2013}.

Next, we prove Theorem \ref{rel-mod8} by using the $2$-adic expansion \eqref{L-2-1}. Recall that $c_k(n)$ in \eqref{L-2-1} denotes the
number of representations of $n$ as a sum of $k$ squares of positive
integers. In particular, $c_1(n)=1$ if $n$ is a square; otherwise, $c_1(n)=0$.

{\noindent \emph{Proof of Theorem \ref{rel-mod8}}}.
 It follows from \eqref{L-2-1} that
\begin{align}\label{P-T-4-1-1}
\overline{p}(n) &\equiv(-1)^n\left(-2c_1(n)+4c_2(n)\right) \pmod{8},
\end{align}
where $n\geq 1$.
Replacing $n$ by $4n$ in \eqref{P-T-4-1-1}, we get
\begin{equation}\label{eq-p4n}
\overline{p}(4n)\equiv-2c_1(4n)+4c_2(4n) \pmod{8}.
\end{equation}
Since $c_1(n)=c_1(4n)$ and $c_2(n)=c_2(4n)$,  \eqref{eq-p4n} 
can be rewritten as 
\begin{align}\label{P-T-4-1-2}
\overline{p}(4n)\equiv-2c_1(n)+4c_2(n) \pmod{8}.
\end{align}
Substituting \eqref{P-T-4-1-2} into \eqref{P-T-4-1-1}, we arrive at
\begin{align*}
\overline{p}(n)\equiv (-1)^{n}\overline{p}(4n)\pmod{8},
\end{align*}
as claimed. \qed

It is easy to see that Corollary \ref{gen-HS} can be obtained by 
iteratively applying Theorem  \ref{rel-mod5} and Theorem \ref{rel-mod8} to  the congruences   $\overline{p}(40n + 35) \equiv 0 \pmod{5}$  and $\overline{p}(40n + 35) \equiv 0 \pmod{8}$ that can be deduced from congruence \eqref{4n+3} by replacing $n$ with $10n+8$.

\section{Proof of Theorem \ref{thm-2}}

In this section, we  prove Theorem \ref{thm-2} by using the Hecke operator on $\phi(q)^3$  along with the fact that $\phi(q)^3$ is a Hecke eigenform in   $M_{\frac{3}{2}}(\tilde{\Gamma}_0(4))$. 

 In view of Theorem \ref{rel-mod5}, to prove Theorem \ref{thm-2}, it suffices to consider the special case $k=0$ that takes the following form.

\begin{thm}
Let $\ell$ be a prime with $\ell\equiv3 \pmod{5}$. Then
\begin{align}\label{l-mod5}
\overline{p}(5\ell^2n)\equiv 0  \pmod{5}
\end{align}
holds for any nonnegative integer $n$ with $\left(\frac{-n}{\ell}\right)=-1$.
\end{thm}
\pf Recall that $\phi(-q)^3$ is a modular form in $ M_{\frac{3}{2}}(\tilde{\Gamma}_0(16))$. Suppose that
\begin{align}\label{eq-exphi}
\phi(-q)^3=\sum_{n\geq0}a(n)q^n
\end{align}
is the Fourier expansion of $\phi(-q)^3$.

Applying the Hecke operator $T_{\frac{3}{2},16}(\ell^2)$ to $\phi(-q)^3$
and using  \eqref{hec-1}, we find that
\begin{equation}\label{heck-1}
\phi(-q)^3\mid T_{\frac{3}{2},16}(\ell^2)=\sum_{n=0}^\infty \left(a(\ell^2n)+\Big(\frac{-n}{\ell}\Big)a(n)+\ell a\Big(\frac{n}{\ell^2}\Big)\right) q^n,
\end{equation}
where $\ell$ is an odd prime.
Replacing $q$ by $-q$ in \eqref{eigenform-1}, 
we see that $\phi(-q)^3$ is a Hecke eigenform in the space $M_{\frac{3}{2}}(\tilde{\Gamma}_0(16))$,  and hence
\begin{align}\label{L-1-2}
\phi(-q)^3\mid T_{\frac{3}{2},16}(\ell^2)=(\ell+1)\phi(-q)^3.
\end{align}
Comparing the coefficients of $q^n$ in \eqref{heck-1} and \eqref{L-1-2}, we deduce that
\begin{align}\label{P-T-5-2}
a(\ell^2n)+\Big(\frac{-n}{\ell}\Big)a(n)+\ell a\Big(\frac{n}{\ell^2}\Big)=(\ell+1)a(n).
\end{align}
Revoking  the  congruence \eqref{L-1-1}, that is,
\begin{equation}\label{eq-phiq3}
\phi(-q)^3\equiv\sum_{n\geq0}\overline{p}(5 n)q^n \pmod{5},
\end{equation}
and comparing \eqref{eq-exphi} with \eqref{eq-phiq3}, we get
\begin{align}\label{P-T-5-3}
a(n)\equiv\overline{p}(5n) \pmod{5}.
\end{align}
Plugging \eqref{P-T-5-3} into \eqref{P-T-5-2}, we deduce that
\begin{align}\label{overpc}
\overline{p}(5\ell^2n)+\left(\frac{-n}{\ell}\right)\overline{p}(5n)+\ell
\overline{p}\left(\frac{5n}{\ell^2}\right)\equiv(\ell+1)\overline{p}(5n) \pmod{5}.
\end{align}
Since $\ell\equiv3 \pmod{5}$ and $(\frac{-n}{\ell})=-1$, we see that $\ell \nmid 5$ and $\ell \nmid n$, so that $\ell^2 \nmid 5n$ and $\overline{p}\left(\frac{5n}{\ell^2}\right)=0$.
Moreover, we have $\left(\frac{-n}{\ell}\right)\equiv(\ell+1)\equiv -1 \pmod{5}$.
Hence congruence \eqref{overpc} becomes
\begin{align*}
\overline{p}(5\ell^2n)\equiv 0 \pmod{5}.
\end{align*}
This completes the proof. \qed

 We now give some special cases of Theorem \ref{thm-2}. Setting $\ell=3$ and $k=0,1$ in Theorem \ref{thm-2}, respectively, we obtain the following
 congruences for $n\geq0$, 
\begin{align*}
&\overline{p}\big(45(3n+1)\big)\equiv 0  \pmod{5},\\[5pt]
&\overline{p}\big(180(3n+1)\big)\equiv 0  \pmod{5}.
\end{align*}
Setting $\ell=13$, $k=0$ in Theorem \ref{thm-2}, we obtain the following congruences for $n\geq0$,
\begin{align*}
&\overline{p}\big(845 (13n+2)\big)\equiv0 \pmod{5},\\[5pt]
&\overline{p}\big(845 (13n+5)\big)\equiv0 \pmod{5},\\[5pt]
&\overline{p}\big(845 (13n+6)\big)\equiv0 \pmod{5},\\[5pt]
&\overline{p}\big(845 (13n+7)\big)\equiv0 \pmod{5},\\[5pt]
&\overline{p}\big(845 (13n+8)\big)\equiv0 \pmod{5},\\[5pt]
&\overline{p}\big(845 (13n+11)\big)\equiv0 \pmod{5}.
\end{align*}

\section{Proof of Theorem \ref{relation-mod5}}

In this section, we complete the proof of Theorem \ref{relation-mod5} by using  the Hecke operator $T_{\frac{3}{2},16}(\ell^2)$ and the Hecke eigenform $\phi(-q)^3$.

{\noindent \emph{Proof of Theorem \ref{relation-mod5}.}}
Setting $\ell=5$ in the congruence relation \eqref{overpc}, we find that
\begin{align}\label{P-T-5-4}
\overline{p}(5n)\equiv \overline{p}(5^3n)+\left(\frac{n}{5}\right)\overline{p}(5n) \pmod{5}.
\end{align}
By the definition of the Legendre symbol, we see that if $n\equiv0 \pmod{5}$, then $\left(\frac{n}{5}\right)=0$.
Hence, by replacing $n$ with $5n$ in   congruence \eqref{P-T-5-4}, we obtain that
\begin{align}
\overline{p}(5^2 n)\equiv\overline{p}(5^4 n) \pmod{5}, \label{P-T-5-6-1}
\end{align}
as claimed. \qed

Furthermore, we note that if $n\equiv\pm1 \pmod{5}$, then
$\left(\frac{n}{5}\right)=1$. Hence by setting $n$ to $5n\pm1$ in \eqref{P-T-5-4}, we deduce that
\begin{align}\label{P-T-5-6-2}
\overline{p}(5^3(5n\pm1))\equiv 0 \pmod{5}.
\end{align}
By iteratively applying the congruence $\overline{p}(5n)\equiv (-1)^{n}\overline{p}(4\cdot 5n) \pmod{5}$ given in Theorem \ref{rel-mod5} and congruence \eqref{P-T-5-6-1} to   \eqref{P-T-5-6-2}, we obtain that
\begin{equation}\label{eq-p45}
\overline{p}(4^k 5^{2i+3}(5n\pm1))\equiv 0 \pmod{5},
\end{equation}
where $n, k,i \ge 0$. This proves
Corollary \ref{coro-2}.

For $n\geq0$, setting $i=0$ and $k=0,1$ in \eqref{eq-p45},  we obtain the following
 special cases 
\begin{align*}
&\overline{p}\big(125(5n\pm 1)\big) \equiv 0 \pmod{5},\\[5pt]
&\overline{p}\big(500(5n\pm 1)\big) \equiv0 \pmod{5}.
\end{align*}

By replacing $n$ by $5n\pm2$ in \eqref{P-T-5-4} and iteratively using the congruence relation \eqref{P-T-5-6-1}, we obtain the following   relation.

\begin{coro}
For $n,i\geq0$, we have
\begin{align*}
\overline{p}\big(5(5n\pm2)\big)\equiv 3\,\overline{p}\big(5^{2i+3}(5n\pm2)\big) \pmod{5}.
\end{align*}
\end{coro}

\vspace{.3cm}

{\noindent \bf Acknowledgments.} This work was supported by the 973 Project, the PCSIRT Project
of the Ministry of Education and the National Science Foundation of China.

\end{document}